\expandafter\ifx\csname mdqon\endcsname\relax
\else \endinput \fi
\message{Document Style Option `german'  Version 2 as of 16 May 1988}
%
\ifx\protect\undefined
\let\protect=\relax \fi
{\catcode`\@=11 
\gdef\allowhyphens{\penalty\@M \hskip\z@skip}

\newcount\U@C\newbox\U@B\newdimen\U@D
\gdef\umlauthigh{\def\"##1{{\accent127 ##1}}}
\gdef\umlautlow{\def\"{\protect\newumlaut}}
\gdef\newumlaut#1{\leavevmode\allowhyphens
     \vbox{\baselineskip\z@skip \lineskip.25ex
     \ialign{##\crcr\hidewidth
     \setbox\U@B\hbox{#1}\U@D .01\p@\U@C\U@D
     \U@D\ht\U@B\advance\U@D -1ex\divide\U@D \U@C
     \U@C\U@D\U@D\the\fontdimen1\the\font
     \multiply\U@D \U@C\divide\U@D 100\kern\U@D
     \vbox to .20ex  
     {\hbox{\char127}\vss}\hidewidth\crcr#1\crcr}}\allowhyphens}
\gdef\highumlaut#1{\leavevmode\allowhyphens
     \accent127 #1\allowhyphens}
\gdef\set@low@box#1{\setbox1\hbox{,}\setbox\z@\hbox{#1}\dimen\z@\ht\z@
     \advance\dimen\z@ -\ht1
     \setbox\z@\hbox{\lower\dimen\z@ \box\z@}\ht\z@\ht1 \dp\z@\dp1}
%
\gdef\@glqq{{\ifhmode \edef\@SF{\spacefactor\the\spacefactor}\else
     \let\@SF\empty \fi \leavevmode
     \set@low@box{''}\box\z@\kern-.04em\allowhyphens\@SF\relax}}
\gdef\glqq{\protect\@glqq}
\gdef\@grqq{\ifhmode \edef\@SF{\spacefactor\the\spacefactor}\else
     \let\@SF\empty \fi \kern-.07em``\kern.07em\@SF\relax}
\gdef\grqq{\protect\@grqq}
\gdef\@glq{{\ifhmode \edef\@SF{\spacefactor\the\spacefactor}\else
     \let\@SF\empty \fi \leavevmode
     \set@low@box{'}\box\z@\kern-.04em\allowhyphens\@SF\relax}}
\gdef\glq{\protect\@glq}
\gdef\@grq{\kern-.07em`\kern.07em}
\gdef\grq{\protect\@grq}
\gdef\@flqq{\ifhmode \edef\@SF{\spacefactor\the\spacefactor}\else
     \let\@SF\empty \fi
     \ifmmode \ll \else \leavevmode
     \raise .2ex \hbox{$\scriptscriptstyle \ll $}\fi \@SF\relax}
\gdef\flqq{\protect\@flqq}
\gdef\@frqq{\ifhmode \edef\@SF{\spacefactor\the\spacefactor}\else
     \let\@SF\empty \fi
     \ifmmode \gg \else \leavevmode
     \raise .2ex \hbox{$\scriptscriptstyle \gg $}\fi \@SF\relax}
\gdef\frqq{\protect\@frqq}
\gdef\@flq{\ifhmode \edef\@SF{\spacefactor\the\spacefactor}\else
     \let\@SF\empty \fi
     \ifmmode < \else \leavevmode
     \raise .2ex \hbox{$\scriptscriptstyle < $}\fi \@SF\relax}
\gdef\flq{\protect\@flq}
\gdef\@frq{\ifhmode \edef\@SF{\spacefactor\the\spacefactor}\else
     \let\@SF\empty \fi
     \ifmmode > \else \leavevmode
     \raise .2ex \hbox{$\scriptscriptstyle > $}\fi \@SF\relax}
\gdef\frq{\protect\@frq}

\global\let\original@ss=\ss
\gdef\newss{\leavevmode\allowhyphens\original@ss\allowhyphens}
\global\let\ss=\newss
\global\let\original@three=\3 
%
%
\gdef\german@dospecials{\do\ \do\\\do\{\do\}\do\$\do\&%
  \do\#\do\^\do\^^K\do\_\do\^^A\do\%\do\~\do\"}
\gdef\german@sanitize{\@makeother\ \@makeother\\\@makeother\$\@makeother\&%
\@makeother\#\@makeother\^\@makeother\^^K\@makeother\_\@makeother\^^A%
\@makeother\%\@makeother\~\@makeother\"}
\global\let\original@dospecials\dospecials
\global\let\dospecials\german@dospecials
\global\let\original@sanitize\@sanitize
\global\let\@sanitize\german@sanitize
\gdef\mdqon{\let\dospecials\german@dospecials
        \let\@sanitize\german@sanitize\catcode`\"\active}
\gdef\mdqoff{\catcode`\"12\let\original@dospecials\dospecials
        \let\@sanitize\original@sanitize}

{\mdqoff
\gdef\@UMLAUT{\"}
\gdef\@MATHUMLAUT{\mathaccent"707F}
\gdef\@SS{\mathchar"7019}
\gdef\dq{"}
}
{\mdqon
\gdef"#1{\if\string#1`\glqq{}%
\else\if\string#1'\grqq{}%
\else\if\string#1a\ifmmode\@MATHUMLAUT a\else\@UMLAUT a\fi
\else\if\string#1o\ifmmode\@MATHUMLAUT o\else\@UMLAUT o\fi
\else\if\string#1u\ifmmode\@MATHUMLAUT u\else\@UMLAUT u\fi
\else\if\string#1A\ifmmode\@MATHUMLAUT A\else\@UMLAUT A\fi
\else\if\string#1O\ifmmode\@MATHUMLAUT O\else\@UMLAUT O\fi
\else\if\string#1U\ifmmode\@MATHUMLAUT U\else\@UMLAUT U\fi
\else\if\string#1e\ifmmode\@MATHUMLAUT e\else\protect \highumlaut e\fi
\else\if\string#1i\ifmmode\@MATHUMLAUT i\else\protect\highumlaut\i \fi
\else\if\string#1E\ifmmode\@MATHUMLAUT E\else\protect\highumlaut E\fi
\else\if\string#1I\ifmmode\@MATHUMLAUT I\else\protect\highumlaut I\fi
\else\if\string#1s\ifmmode\@SS\else\ss\fi{}%
\else\if\string#1-\allowhyphens\-\allowhyphens
\else\if\string#1\string"\hskip\z@skip
\else\if\string#1|\discretionary{-}{}{\kern.03em}%
\else\if\string#1c\allowhyphens\discretionary{k-}{}{c}\allowhyphens
\else\if\string#1f\allowhyphens\discretionary{ff-}{}{f}\allowhyphens
\else\if\string#1l\allowhyphens\discretionary{ll-}{}{l}\allowhyphens
\else\if\string#1m\allowhyphens\discretionary{mm-}{}{m}\allowhyphens
\else\if\string#1n\allowhyphens\discretionary{nn-}{}{n}\allowhyphens
\else\if\string#1p\allowhyphens\discretionary{pp-}{}{p}\allowhyphens
\else\if\string#1t\allowhyphens\discretionary{tt-}{}{t}\allowhyphens
\else\if\string#1<\flqq{}%
\else\if\string#1>\frqq{}%
\else             \dq #1%
\fi\fi\fi\fi\fi\fi\fi\fi\fi\fi\fi\fi\fi\fi\fi\fi\fi\fi\fi\fi\fi\fi\fi\fi\fi}
}
\gdef\dateaustrian{\def\today{\number\day.~\ifcase\month\or
  J\"anner\or Februar\or M\"arz\or April\or Mai\or Juni\or
  Juli\or August\or September\or Oktober\or November\or Dezember\fi
  \space\number\year}}
\gdef\dategerman{\def\today{\number\day.~\ifcase\month\or
  Januar\or Februar\or M\"arz\or April\or Mai\or Juni\or
  Juli\or August\or September\or Oktober\or November\or Dezember\fi
  \space\number\year}}
\gdef\dateUSenglish{\def\today{\ifcase\month\or
 January\or February\or March\or April\or May\or June\or
 July\or August\or September\or October\or November\or December\fi
 \space\number\day, \number\year}}
\gdef\dateenglish{\def\today{\ifcase\day\or
 1st\or 2nd\or 3rd\or 4th\or 5th\or
 6th\or 7th\or 8th\or 9th\or 10th\or
 11th\or 12th\or 13th\or 14th\or 15th\or
 16th\or 17th\or 18th\or 19th\or 20th\or
 21st\or 22nd\or 23rd\or 24th\or 25th\or
 26th\or 27th\or 28th\or 29th\or 30th\or
 31st\fi
 ~\ifcase\month\or
 January\or February\or March\or April\or May\or June\or
 July\or August\or September\or October\or November\or December\fi
 \space \number\year}}
\gdef\datefrench{\def\today{\ifnum\day=1\relax 1\/$^{\rm er}$\else
  \number\day\fi \space\ifcase\month\or
  janvier\or f\'evrier\or mars\or avril\or mai\or juin\or
  juillet\or ao\^ut\or septembre\or octobre\or novembre\or d\'ecembre\fi
  \space\number\year}}
%
%
\gdef\captionsgerman{%
\def\refname{Literatur}%
\def\abstractname{Zusammenfassung}%
\def\bibname{Literaturverzeichnis}%
\def\chaptername{Kapitel}%
\def\appendixname{Anhang}%
\def\contentsname{Inhaltsverzeichnis}%
\def\listfigurename{Abbildungsverzeichnis}%
\def\listtablename{Tabellenverzeichnis}%
\def\indexname{Index}%
\def\figurename{Abbildung}%
\def\tablename{Tabelle}%
\def\partname{Teil}}
\gdef\captionsenglish{%
\def\refname{References}%
\def\abstractname{Abstract}%
\def\bibname{Bibliography}%
\def\chaptername{Chapter}%
\def\appendixname{Appendix}%
\def\contentsname{Contents}%
\def\listfigurename{List of Figures}%
\def\listtablename{List of Tables}%
\def\indexname{Index}%
\def\figurename{Figure}%
\def\tablename{Table}%
\def\partname{Part}}
\gdef\captionsfrench{%
\def\refname{R\'ef\'erences}%
\def\abstractname{R\'esum\'e}%
\def\bibname{Bibliographie}%
\def\chaptername{Chapitre}%
\def\appendixname{Appendice}%
\def\contentsname{Table des mati\`eres}%
\def\listfigurename{Liste des figures}%
\def\listtablename{Liste des tables}%
\def\indexname{Index}%
\def\figurename{Figure}%
\def\tablename{Table}%
\def\partname{Partie}}%
\newcount\language 
\newcount\USenglish  \global\USenglish=0
\newcount\german     \global\german=1
\newcount\austrian   \global\austrian=2
\newcount\french     \global\french=3
\newcount\english    \global\english=4
\gdef\setlanguage#1{\language #1\relax
  \expandafter\ifcase #1\relax
  \dateUSenglish  \captionsenglish   \or
  \dategerman     \captionsgerman    \or
  \dateaustrian   \captionsgerman    \or
  \datefrench     \captionsfrench    \or
  \dateenglish    \captionsenglish   \fi}
\gdef\originalTeX{\mdqoff \umlauthigh
  \let\ss\original@ss \let\3\original@three
  \setlanguage{\USenglish}}
\gdef\germanTeX{\mdqon \umlautlow \let\ss\newss \let\3\ss
  \let\dospecials\german@dospecials
  \setlanguage{\german}}
} 
%
%
\germanTeX
%
%

\output={\if N\header\headline={\hfill}\fi
\plainoutput\global\let\header=Y}
\magnification\magstep1
\tolerance = 500
\hsize=14.4true cm
\vsize=22.5true cm
\parindent=6true mm\overfullrule=2pt
\newcount\kapnum \kapnum=0
\newcount\parnum \parnum=0
\newcount\procnum \procnum=0
\newcount\nicknum \nicknum=1
\font\ninett=cmtt9

\font\ninebf=cmbx9

\font\sixbf=cmbx6
\font\ninesl=cmsl9

\font\nineit=cmti9

\font\ninerm=cmr9

\font\sixrm=cmr6
\font\ninei=cmmi9
\font\eighti=cmmi8
\font\sixi=cmmi6
\skewchar\ninei='177 \skewchar\eighti='177 \skewchar\sixi='177
\font\ninesy=cmsy9
\font\eightsy=cmsy8
\font\sixsy=cmsy6
\skewchar\ninesy='60 \skewchar\eightsy='60 \skewchar\sixsy='60
\font\titelfont=cmr10 scaled 1440
\font\paragratit=cmbx10 scaled 1200

\font\name=cmcsc10
\font\emph=cmbxti10

\font\tenmsbm=msbm10
\font\sevenmsbm=msbm7
%

%

%
\font\teneufm=eufm10
\font\seveneufm=eufm7
\font\fiveeufm=eufm5
\newfam\eufmfam
\textfont\eufmfam=\teneufm
\scriptfont\eufmfam=\seveneufm
\scriptscriptfont\eufmfam=\fiveeufm

\font\tenmsam=msam10
\font\sevenmsam=msam7
\font\fivemsam=msam5
\newfam\msamfam
\textfont\msamfam=\tenmsam
\scriptfont\msamfam=\sevenmsam
\scriptscriptfont\msamfam=\fivemsam
\font\tenmsbm=msbm10
\font\sevenmsbm=msbm7
\font\fivemsbm=msbm5
\newfam\msbmfam
\textfont\msbmfam=\tenmsbm
\scriptfont\msbmfam=\sevenmsbm
\scriptscriptfont\msbmfam=\fivemsbm
\def\Bbb#1{{\fam\msbmfam\relax#1}}
\def\cz{{\kern0.4pt\Bbb C\kern0.7pt}
}
\def\ez{{\kern0.4pt\Bbb E\kern0.7pt}
}
\def\fz{{\kern0.4pt\Bbb F\kern0.3pt}}
\def\gz{{\kern0.4pt\Bbb Z\kern0.7pt}}
\def\hz{{\kern0.4pt\Bbb H\kern0.7pt}
}
\def\kz{{\kern0.4pt\Bbb K\kern0.7pt}
}
\def\nz{{\kern0.4pt\Bbb N\kern0.7pt}
}
\def\oz{{\kern0.4pt\Bbb O\kern0.7pt}
}
\def\rz{{\kern0.4pt\Bbb R\kern0.7pt}
}
\def\sz{{\kern0.4pt\Bbb S\kern0.7pt}
}
\def\pz{{\kern0.4pt\Bbb P\kern0.7pt}
}
\def\qz{{\kern0.4pt\Bbb Q\kern0.7pt}
}
\newskip\ttglue
\def\ninepoint{\def\rm{\fam0\ninerm}%
  \textfont0=\ninerm \scriptfont0=\sixrm \scriptscriptfont0=\fiverm
  \textfont1=\ninei \scriptfont1=\sixi \scriptscriptfont1=\fivei
  \textfont2=\ninesy \scriptfont2=\sixsy \scriptscriptfont2=\fivesy
  \textfont3=\tenex \scriptfont3=\tenex \scriptscriptfont3=\tenex
  \def\it{\fam\itfam\nineit}%
  \textfont\itfam=\nineit
  \def\sl{\fam\slfam\ninesl}%
  \textfont\slfam=\ninesl
  \def\bf{\fam\bffam\ninebf}%
  \textfont\bffam=\ninebf \scriptfont\bffam=\sixbf
   \scriptscriptfont\bffam=\fivebf
  \def\tt{\fam\ttfam\ninett}%
  \textfont\ttfam=\ninett
  \tt \ttglue=.5em plus.25em minus.15em
  \normalbaselineskip=11pt
  \font\name=cmcsc9
  \let\sc=\sevenrm
  \let\big=\ninebig
  \setbox\strutbox=\hbox{\vrule height8pt depth3pt width0pt}%
  \normalbaselines\rm
  \def\sl{\it}}

\headline={\ifodd\pageno\rightheadline\else\leftheadline\fi}
\def\rightheadline{\ninepoint Paragraphen"uberschrift\hfill\folio}
\def\leftheadline{\ninepoint\folio\hfill Chapter"uberschrift}
\let\header=Y
\def\titel#1{\need 9cm \vskip 2truecm
\parnum=0\global\advance \kapnum by 1
{\baselineskip=16pt\lineskip=16pt\rightskip0pt
plus4em\spaceskip.3333em\xspaceskip.5em\pretolerance=10000\noindent
\titelfont Chapter \uppercase\expandafter{\romannumeral\kapnum}.
#1\vskip2true cm}\def\leftheadline{\ninepoint
\folio\hfill Chapter \uppercase\expandafter{\romannumeral\kapnum}.
#1}\let\header=N
}
\def\Titel#1{\need 9cm \vskip 2truecm
\global\advance \kapnum by 1
{\baselineskip=16pt\lineskip=16pt\rightskip0pt
plus4em\spaceskip.3333em\xspaceskip.5em\pretolerance=10000\noindent
\titelfont\uppercase\expandafter{\romannumeral\kapnum}.
#1\vskip2true cm}\def\leftheadline{\ninepoint
\folio\hfill\uppercase\expandafter{\romannumeral\kapnum}.
#1}\let\header=N
}
\def\need#1cm {\par\dimen0=\pagetotal\ifdim\dimen0<\vsize
\global\advance\dimen0by#1 true cm
\ifdim\dimen0>\vsize\vfil\eject\noindent\fi\fi}
\def\neupara#1{\par\penalty-2000
\procnum=0\global\advance\parnum by 1
\vskip1cm\noindent{\paragratit \the\parnum. #1}%
\def\rightheadline{\ninepoint\S\the\parnum.\ #1\hfill \folio}%
\vskip 8mm\noindent}
\def\Proclaim #1 #2\finishproclaim {\bigbreak\noindent
{\bf#1\unskip{}. }{\it#2}\medbreak\noindent}
%
\gdef\proclaim #1 #2 #3\finishproclaim {\bigbreak\noindent%
\global\advance\procnum by 1
{%
{\relax\ifodd \nicknum
\hbox to 0pt{\vrule depth 0pt height0pt width\hsize
   \quad \ninett#3\hss}\else {}\fi}%
\bf\the\parnum.\the\procnum\ #1\unskip{}. }
{\it#2}
\immediate\write\num{\string\def
 \expandafter\string\csname#3\endcsname
 {\the\parnum.\the\procnum}}
\medbreak\noindent}
\newcount\stunde \newcount\minute \newcount\hilfsvar
\def\uhrzeit{
    \stunde=\the\time \divide \stunde by 60
    \minute=\the\time
    \hilfsvar=\stunde \multiply \hilfsvar by 60
    \advance \minute by -\hilfsvar
    \ifnum\the\stunde<10
    \ifnum\the\minute<10
    0\the\stunde:0\the\minute~Uhr
    \else
    0\the\stunde:\the\minute~Uhr
    \fi
    \else
    \ifnum\the\minute<10
    \the\stunde:0\the\minute~Uhr
    \else
    \the\stunde:\the\minute~Uhr
    \fi
    \fi
    }

 \def\calF{{\cal F}}

\def\calK{{\cal K}} 
\def\calM{{\cal M}} \def\calN{{\cal N}}

 \def\calZ{{\cal Z}}

\def\GL{\mathop{\rm GL}\nolimits}

\def\id{\mathop{\rm id}\nolimits}
\def\im{\mathop{\rm Im}\nolimits} \def\Im{\im}

\def\kernel{\mathop{\rm kernel}\nolimits}

\def\mod{\mathop{\rm mod}\nolimits}

\def\Sp{\mathop{\rm Sp}\nolimits}

\def\Sym{\mathop{\rm Sym}\nolimits}
\def\boxit#1{
  \vbox{\hrule\hbox{\vrule\kern6pt
  \vbox{\kern8pt#1\kern8pt}\kern6pt\vrule}\hrule}}
\def\Boxit#1{
  \vbox{\hrule\hbox{\vrule\kern2pt
  \vbox{\kern2pt#1\kern2pt}\kern2pt\vrule}\hrule}}

\def\smallni{\smallskip\noindent }
\def\medni{\medskip\noindent }

\def\lo{\longrightarrow}

\def\loma{\longmapsto}

\def\pii{\pi {\rm i}}

\def\square{\hbox{\hbox to 0pt{$\sqcup$\hss}\hbox{$\sqcap$}}}
\def\qed{\ifmmode\square\else{\unskip\nobreak\hfil
\penalty50\hskip3em\null\nobreak\hfil\square
\parfillskip=0pt\finalhyphendemerits=0\endgraf}\fi}
\def\pn{\the\parnum.\the\procnum}
\def\downmapsto{{\buildrel
        {\vbox{\hbox{\hskip.2pt$\scriptstyle-$}}}
        \over{\raise7pt\vbox{\vskip-4pt\hbox{$\textstyle\downarrow$}}}}}
\nopagenumbers
\font\pro =cmss10
\immediate\newwrite\num
\let\header=N
\def\transpose#1{\kern1pt{^t\kern-1pt#1}}%
\input standard.num
\immediate\openout\num=standard.num
\immediate\newwrite\num\immediate\openout\num=standard.num
\def\RAND#1{\vskip0pt\hbox to 0mm{\hss\vtop to 0pt{%
  \raggedright\ninepoint\parindent=0pt%
  \baselineskip=1pt\hsize=2cm #1\vss}}\noindent}
\noindent
in memoriam\hfill\break
Jun--Ichi Igusa (1924--2013)
\vskip5mm\noindent
\centerline{\titelfont Basic vector valued Siegel modular forms of genus two}%
\def\leftheadline{\ninepoint\folio\hfill
Vector valued modular forms}%
\def\rightheadline{\ninepoint Introduction\hfill \folio}%
\headline={\ifodd\pageno\rightheadline\else\leftheadline\fi}
\vskip 1.5cm
\leftline{\it \hbox to 6cm{Eberhard Freitag\hss}
Riccardo Salvati
Manni  }
  \leftline {\it  \hbox to 6cm{Mathematisches Institut\hss}
Dipartimento di Matematica, }
\leftline {\it  \hbox to 6cm{Im Neuenheimer Feld 288\hss}
Piazzale Aldo Moro, 2}
\leftline {\it  \hbox to 6cm{D69120 Heidelberg\hss}
 I-00185 Roma, Italy. }
\leftline {\tt \hbox to 6cm{freitag@mathi.uni-heidelberg.de\hss}
salvati@mat.uniroma1.it}
\vskip1cm
\centerline{\paragratit \rm  2013}%
\vskip5mm\noindent%
\let\header=N%
\def\St{{\rm St}}%
\def\grad{{\rm grad}}%
\vskip1.5cm\noindent
{\paragratit Introduction}%
\medni
For a congruence subgroup $\Gamma\subset\Sp(n,\gz)$ of the Siegel modular group of genus
$n$, a starting weight $r_0>0$, $2r_0\in\gz$,  and a starting multiplier system $v$
of weight $r_0$
one can consider the ring of modular forms (for definitions we refer to Sect.~1)
$$A(\Gamma)=A(\Gamma,r_0,r)=\bigoplus_{r\in\gz}[\Gamma,rr_0,v^r].$$
This is a finitely generated algebra. In addition, let $\varrho:\GL(n,\cz)\to\GL(\calZ)$
be a rational representation on a finite dimensional complex vector space $\calZ$.
We always will assume that $\varrho$ is irreducible and polynomial and  does not
vanish on the determinant surface $\det A=0$.
Then we consider the space
$[\Gamma,rr_0,v^r,\varrho]$
of all vector-valued holomorphic modular forms $f:\hz_n\to\calZ$ of transformation type
$$f(MZ)=v^r(M)\det(CZ+D)^{r_0r}\varrho(CZ+D)f(Z).$$
We collect these spaces to the graded $A(\Gamma)$-module
$$\calM=\calM_\Gamma(r_0,v,\varrho):=\bigoplus_{r\in\gz}[\Gamma,rr_0,v^r,\varrho].$$
There are twisted variants of these modules [Wi]. For a character $\chi$ on $\Gamma$
one can consider
$$\calM^\chi=\bigoplus_{r\in\gz}[\Gamma,rr_0,\chi v^r,\varrho].$$
These modules are finitely generated. It is  a natural task to look for examples
where the structure of this module can be determined.
\smallskip
Meanwhile there appeared several papers
getting results into this direction using different methods,
[Ao, CG, Do, Ib, Sa, Sat, Wi]. 
Our method is a further development of Wieber's geometric method [Wi] which he used
to solve certain $\Gamma_2[2,4]$-cases. 
\smallskip
The geometric method of Wieber rests on the fact that
vector valued modular forms sometimes can be interpreted as $\Gamma$-invariant tensors
and hence as (usually rational) tensors on the Siegel modular variety
$\hz_n/\Gamma$. In some cases the structure of this variety is known which enables to study
tensors on it in detail. 
A vector valued modular form can define a tensor only if the representation
$\varrho$ up to a power of the determinant occurs in some tensor power of the representation
$\Sym^2$. This is not always the case.
For example the standard representation
of $\GL(2,\cz)$ does not have this property. In this paper we describe
a modification of Wieber's method which allows to recover his main results in [Wi]
in a quick way
and which applies to more cases as for example the standard representation. 
\smallskip
Recall that the principal congruence  subgroup is defined as
$$\Gamma_n[q]=\kernel\bigl(\Sp(n,\gz)\lo\Sp(n,\gz/q\gz)\bigr)$$
and Igusa's subgroup as
$$\Gamma_n[q,2q]:=\left\{M\in\Gamma_n[q],\quad(C\transpose D)_0
\equiv(A\transpose B)_0\equiv0\;\mod2q\;\right\}.$$
Here $S_0$ denotes the column built of the diagonal of a square matrix $S$.
\smallskip
Besides Wieber's known results we treat
in this paper a new example that belongs to the group
$\Gamma_2[4,8]$. The starting weight is $1/2$, the starting multiplier system is
the theta multiplier system $v_\vartheta$ and for $\varrho$ we take the standard representation.
In this case we will determine the structure of $\calM$ completely (Theorem \MainT).
It will turn out that $\calM$  can be generated by the $\Gamma_2$ orbits of two specific
modular forms. We will describe the relations and -- as a consequence -- we will obtain
the Hilbert function of this module (Theorem \MainT).
\smallskip
We want to thank Wieber for fruitful discussion and for his help with quite
involved computer calculations.
\neupara{Vector valued modular forms}%
We consider the Siegel modular group $\Gamma_n=\Sp(n,\gz)$ of genus $n$. It consists
of all integral $2n\times 2n$-matrices $M$ such that $\transpose MIM=I$,
where $I={0\,-E\choose E\,\phantom{-}0}$ ($E$ denotes the unit matrix) 
is the standard alternating matrix.
Let $\Gamma\subset\Sp(n,\gz)$ be a congruence subgroup,
$r$ be an integer, $v$ a multiplier system of weight $r/2$ on $\Gamma$ and
$\varrho:\GL(n,\cz)\to\GL(n,\calZ)$ a rational representation on a finite dimensional
complex vector space $\calZ$. We assume that $\varrho$ is reduced which means that it is
polynomial and does not vanish along the determinant surface $\det(A)=0$.
Then we can consider the space $[\Gamma,r/2,v,\varrho]$. It consists
of holomorphic functions $f:\hz_n\to \calZ$ on the Siegel upper half-plane
$$\hz_n=\bigl\{\,Z\in\cz^{n\times n};\quad \Im Z>0\ \hbox{(positive definit)}\,\bigr\}$$
with the transformation property
$$f(MZ)=v(M)\sqrt{\det(CZ+D)}^r\varrho(CZ+D)f(Z)\qquad(M\in\Gamma).$$
(In the case $n=1$ a condition at the cusps has to be added.)
We can also consider meromorphic solutions $f$ and call them meromorphic modular forms
if they satisfy a meromorphicity condition at the cusps (which in most cases will be
automatically true). What we demand is that there exists a non-vanishing holomorphic
scalar valued form $g$ that that $fg$ is holomorphic.
We denote the space of meromorphic modular forms by
$\{\Gamma,r/2,v,\varrho\}$. It is a vector space of dimension $\le{{\rm Rank}}(\varrho)$
over the field
of modular functions
$$K(\Gamma)=\{\Gamma,0,\hbox{triv},\hbox{triv}\}.$$
In the case that $\varrho$ is the one-dimensional trivial representation,
we simply write $\{\Gamma,r/2,v\}$ instead of $\{\Gamma,r/2,v,\varrho\}$
and similarly we skip $v$ if $r$ is even and $v$ trivial. The same convention is
used for the spaces of holomorphic modular forms.
\neupara{Thetanullwerte}%
An element $m\in\{0,1\}^{2n}$ is called a theta characteristic of genus $n$.
Usually it is considered as column and divided into to columns $a,b$ of length $n$.
It is called even if $\transpose ab$ is even and odd else. We use the classical
theta series
$$\vartheta[m](Z,z)=\sum_{g\in\gz^n}\exp\pii(Z[g+a/2]+2\transpose(g+a/2)(z+b/2).$$
We are interested in the nullwerte
$$\vartheta[m](Z)=\vartheta[m](Z,0)$$ 
and in the nullwerte of the derivatives
$${\partial\over\partial z_i}\vartheta[m](Z,z)\big\vert_{z=0}.$$
We collect them in a column which we denote by $\grad\vartheta(Z)$.
We recall that the theta nullwerte are non-zero only for even  and
the gradients for odd characteristics. 
\smallskip
Besides the nullwerte of first kind $\vartheta[m]$ the nullwerte of second kind
$$f_a(Z):=\vartheta\Bigl[{a\atop 0}\Bigr](2Z),\qquad a\in(\gz/2\gz)^n,$$
will play a role.
\smallskip
We recall that $\vartheta[0](Z)$ is a modular form of weight $1/2$ for the
{\it theta group}
$$\Gamma_{n,\vartheta}:=\Gamma_n[1,2]$$
with respect to a certain multiplier system $v_\vartheta$ on this group.
Since $\Gamma_{n,\vartheta}\supset\Gamma_n[2]$ we have in particular
$$\vartheta[0]\in[\Gamma_n[2],1/2,v_\vartheta].$$
For each characteristic there
exists a character $\chi_m$ on $\Gamma_n[2]$ which is trivial on $\Gamma_n[4,8]$
and quadratic on the group $\Gamma_n[2,4]$.
such that
$$\vartheta[m]\in[\Gamma_n[2],1/2,v_\vartheta\chi_m]\quad
\hbox{and}\quad
\grad\vartheta[m]\in[\Gamma_n[2],1/2,v_\vartheta\chi_m,\St].$$
In particular, all thetanullwerte $\vartheta[m]$ are contained in
$[\Gamma_n[4,8],1/2,v_\vartheta]$.
For details we refer to [SM].
\smallskip
Similar results hold for the thetas of second kind $f_a$. They are modular forms
for $\Gamma_n[2,4]$ with respect to a certain multiplier system $v_\Theta$ on this group,
$$f_a\in[\Gamma_n[2,4],1/2,v_\Theta].$$
We consider the rings
$$A(\Gamma_n[4,8])=\bigoplus_{r\in\gz}[\Gamma_n[4,8],r/2,v_\vartheta^r],\quad
A(\Gamma_n[2,4])=\bigoplus_{r\in\gz}[\Gamma_n[2,4],r/2,v_\Theta^r].
$$
So the starting weights are $1/2$ in both cases and the starting multiplier system is
$v_\vartheta$ for $\Gamma_n[4,8]$ but $v_\Theta$ for $\Gamma_n[2,4]$.
We mention that the two multiplier systems are different on $\Gamma_n[4,8]$.
The following results are basic. The first one has proved by Igusa 1964 [Ig] 
the second by Runge 1994 [Ru1,Ru2].
\proclaim
{Theorem  (Igusa, Runge)}
{$$\eqalign{
\hbox to 5cm{$A(\Gamma_n[4,8])=\cz[\dots,\vartheta[m],\dots]$\hfil}&
\quad\hbox{for}\quad n\le 2,\cr
\hbox to 5cm{$A(\Gamma_n[2,4])=\cz[\dots,f_a,\dots]$\hfil}&
\quad\hbox{for}\quad n\le 3.\cr}$$}
IgRu%
\finishproclaim
In the case $n=2$ Runge obtained an even  better result. 
It is known that the square of $v_\Theta$ is a non-trivial quadratic character 
on $\Gamma_2[2,4]$. 
The kernel of $v_\Theta^2$ is a subgroup of index two of $\Gamma_2[2,4]$. We use Runge's
notation $\Gamma_2^*[2,4]$ for it. 
\smallskip 
Recall that Igusa's modular form $\chi_5$  is the unique cusp form of weight
5 for the full Siegel modular form. It can be defined as the product of the ten
theta constants of first kind. Its character is trivial on $\Gamma_2[2]$.
\smallskip
We denote the 4 functions $f_a$ in the ordering $(0,0)$, $(0,1)$, $(1,0)$, $(1,1)$
by $f_0,\dots,f_3$. 
\proclaim
{Theorem (Runge)}
{
$$A(\Gamma^*_2[2,4])=\bigoplus_{r\in\gz}[\Gamma^*_2[2,4],r/2,v_\Theta^r]=
\cz[f_0,\dots,f_3]\oplus\chi_5\cz[f_0,\dots,f_3].$$}
RunZ%
\finishproclaim
\neupara{A first example due to Wieber}%
From now on we shall assume $n=2$.
We consider the module $\calM$ introduced in the introduction, for
the group $\Gamma_2[2,4]$, starting weight $1/2$ and starting multiplier system
$v_\Theta$ and for the representation $\Sym^2$.
We will write $\calM^+$ instead of $\calM$ since, in the next section, we shall
treat a twisted variant $\calM^-$.
\smallskip
The representation $\Sym^2$ of $\GL(2,\cz)$.
can be realized on the space of symmetric $2\times 2$-matrices and the action
of $\GL(2,\cz)$ is given by $AW\transpose A$. 
\smallskip
This means that we have to consider
symmetric $2\times2$-matrices $f$ of holomorphic functions with the
transformation property
$$f(MZ)=v_\Theta(M)^r\sqrt{CZ+D)}^r(CZ+D)f(Z)\transpose(CZ+D).$$
They define a vector space 
$$\calM^+(r)=[\Gamma_2[2,4]),r/2,v_\Theta^r,\Sym^2].$$
We consider the direct
sum
$$\calM^+=\bigoplus_{r\in\gz}\calM^+(r).$$
This is a a graded module over the ring 
$$A(\Gamma_2[2,4]):=\bigoplus[\Gamma,r/2,v_\Theta^r]=\cz[f_0,\dots,f_3].$$
Elements of $\calM^+$ can be constructed as follows. Let 
$f,g\in [\Gamma,r/2,v_\Theta^r]$, $g\ne 0$. Then $f/g$ is a modular function and
$d(f/g)$ is a meromorphic differential. It can be considered as element
of $\{\Gamma_2[2,4],0,\Sym^2\}$. If we multiply it by $g^2$ we get a
holomorphic form
$$\{f,g\}=g^2d(f/g)\in\calM^+(2r).$$
\proclaim
{Theorem (Wieber)}
{We have
$$\calM^+=\sum_{0\le i<j\le 3}\cz[f_0,\dots,f_3]\{f_i,f_j\}.$$
Defining relations of this module are
$$f_k[f_i,f_j]=f_j[f_i,f_k]+f_i[f_k,f_j],\quad [f_i,f_j]+[f_j,f_i]=0.$$
}
WiebT%
\finishproclaim
{\it Proof.\/} We give a new simple proof for this result.
We consider $\{\Gamma[2,4],2,\Sym^2\}$ as vector space
over the field of modular functions. The
three forms $\{f_0,f_i\}$, $1\le i\le 3$, give a basis of this vector space. 
\proclaim
{Lemma}
{With some constant $C$ we have
$$f_0^4d(f_1/f_0)\wedge d(f_2/f_0)\wedge d(f_3/f_0)=C\chi_5dz_0\wedge dz_1\wedge dz_2,\qquad
Z=\pmatrix{z_0&z_1\cr z_1&z_2}.$$
{\bf Corollary.}
If we consider the three elements
$$\{f_0,f_1\},\quad \{f_0,f_2\},\quad \{f_0,f_3\}$$
as a $3\times 3$-matrix, its determinant is up to a constant factor
$f_0^2\chi_5$.}
DetSyme%
\finishproclaim
This lemma is well-known.
For sake of completeness we give the  argument.
The left hand side is holomorphic and it defines
a modular form of weight 5 with respect to the full Siegel modular form. So it must
be a constant multiple of $\chi_5$.\qed
\smallskip
We mentioned already in the introduction that $\{\Gamma,r/2,v,\varrho\}$
is a vector space of dimension $\le{{\rm Rank}}(\varrho)$. We can get a basis of this space
if we multiply the three $\{f_0,f_i\}$ by $f_0^{r-1}$.
Hence an arbitrary $T\in \{\Gamma,r/2,v,\varrho\}$ can be written in the form
$$T=g_1\{f_0,f_1\}+g_2\{f_0,f_2\}+g_3\{f_0,f_3\}$$
where $g_i$ are meromorphic modular forms in $\{\Gamma_2[2,4],(r-2)/2,v_\Theta^r\}$.
They are rational functions in the $f_i$. If $T$ is holomorphic then
$f_0^2\chi_5 g_i$ must be holomorphic. 
The form $f_0^2\chi_5 g_i$ is contained in $A(\Gamma_2^*[2,4])$ 
and $f_0^2g_i$ are also rational in the $f_i$.
From Runge's result Theorem \RunZ\ follows that
$f_0^2g_i$ is holomorphic. In other words
$$\calM^+\subset {1\over f_0^2}\sum_{i=1}^3\cz[f_0,\dots,f_3]\{f_0,f_i\}.$$
The rest is just Wieber's argument that we can permute the variables
and obtain that $\calM$ is contained in the intersection of
4 modules, 
$$\calM^+\subset \bigcap_{i=0}^3{1\over f_i^2}\sum_{j=1}^3\cz[f_0,\dots,f_3]\{f_i,f_j\}.$$
Using the fact that $\calM^+$ is contained in the free module
generated by $df_i$, 
it is easy to show (compare [Wi]) that this intersection equals
$$\sum_{0\le i<j\le 3}\cz[f_0,\dots,f_3]\{f_i,f_j\}.$$
which is Wieber's result.
\neupara{A second example of Wieber}%
Wieber also considers the twist of $\calM^+$ with the quadratic character
$v_\Theta^2$. To be precise, he introduces the spaces $\calM^-(r)$ consisting
of holomorphic forms of the type
$$f(MZ)=v_\Theta(M)^{r+2}\sqrt{CZ+D)}^r(CZ+D)f(Z)\transpose(CZ+D).$$
They can be collected to 
$$\calM^-=\bigoplus\calM^-(r)$$
which is also a graded module over $\cz[f_0,\dots,f_4]$. There are obvious inclusions
$$\chi_5\calM^+\subset\calM^-,\quad \chi_5\calM^-\subset\calM^+.$$
Following more general constructions of Ibukiyama [Ib],
Wieber defined elements of $\calM^-$ in a different way.
He considers three homogenous elements $f,g,h$ of $\cz[f_0,\dots,f_3]$ of degree
$r$ and considers then the differential form
$$d(g/f)\wedge d(h/f)=h_0dz_1\wedge dz_2+h_1 dz_0\wedge dz_2+h_2dz_0\wedge dz_1.$$
Then 
$$\pmatrix{h_2&-h_1\cr -h_1&h_0}\in \{\Gamma_2[2,4], 1, \Sym^2\}.$$
(The multiplier system is trivial.)  We set
$$\{f,g,h\}=f^3\pmatrix{h_2&-h_1\cr -h_1&h_0}.$$
It is easy to see that this is holomorphic. It is contained in 
$$[\Gamma_2[2,4],(3r+2)/2,v_\Theta^{3r},\Sym^2]=\calM^-(3r+2).$$
If $f,g,h$ have weight $1/2$, then this form is contained in $\calM^-(5)$.
%
%
%
An arbitrary $T\in [\Gamma,(r+2)/2,v,\varrho]$ can be written in the form
$$T=g_1\{f_0,f_1,f_2\}+g_2\{f_0,f_1,f_3\}+g_3\{f_0,f_2,f_3\}$$
where $g_i$ are meromorphic modular forms in $\{\Gamma_2[2,4],r/2,v_\Theta^r\}$.
They are expressible as quotients of homogenous polynomials in the variables $f_i$.
We have to work out that this form, or equivalently
$$f_0^3(g_1d(f_1/f_0)\wedge d(f_2/f_0)+g_2d(f_1/f_0)\wedge d(f_3/f_0)+
g_3d(f_2/f_0)\wedge d(f_3/f_0)),$$
is holomorphic. We take the wedge product with $f_0^2df(f_i/f_0)$ and obtain from 
Lemma \DetSyme\ that
$f_0\chi_5 g_i$ are holomorphic. Hence the argument of the previous section shows
that $f_0g_i$ is a polynomial in the $f_i$.
%
%
Similar to the previous section we obtain
$$\calM^-\subset\bigcap_{i=0}^3 {1\over f_i}\sum_{i<j<k}\cz[f_0,\dots,f_3]\{f_i,f_j,f_k\}.$$
A simple argument now gives the following result.
\proclaim
{Theorem (Wieber)}
{We have
$$\calM^-=\sum_{0\le i<j<k\le 3}\cz[f_0,\dots,f_3]\{f_i,f_j,f_k\}.$$
Defining relation of this module is
$$f_3\{f_0,f_1,f_2\}=f_0\{f_1,f_2,f_3\}-f_1\{f_0,f_2,f_3\}+f_2\{f_0,f_1,f_3\}.$$
}
WiebTz%
\finishproclaim
\neupara{The standard representation}%
In this section we study the module $\calM$ for the group $\Gamma_2[4,8]$,
starting weight $1/2$, starting multiplier system $v_\vartheta$ and the standard representation
$\St=\id$,
$$\calM=\bigoplus [\Gamma_2[4,8] r/2 , v_{\vartheta}^r, \St].$$
This is a module over the ring $A(\Gamma_2[4,8])$ which, by Igusa's result, is generated by the
ten even theta nullwerte.
We will order them as follows:
$$(m^{(1)},\dots,m^{(10)})=\pmatrix{
0&0&0&0&0&0&1&1&1&1\cr
0&0&0&0&1&1&0&0&1&1\cr
0&0&1&1&0&1&0&0&0&1\cr
0&1&0&1&0&0&0&1&0&1\cr
}$$
The associated theta series are denote by $\vartheta_1,\dots,\vartheta_{10}$ (in this ordering).
They satisfy the quartic Riemann relations which are defining relations:\smallni
\hbox{
\vbox{\hsize=5.5cm
$$\leqalignno{&
\vartheta_6^2\vartheta_8^2-\vartheta_4^2\vartheta_9^2+\vartheta_1^2\vartheta_{10}^2=0,\cr&
\vartheta_5^2\vartheta_8^2-\vartheta_2^2\vartheta_9^2+\vartheta_3^2\vartheta_{10}^2=0,\cr&
\vartheta_7^4-\vartheta_8^4-\vartheta_9^4+\vartheta_{10}^4=0,\cr&
\vartheta_6^2\vartheta_7^2-\vartheta_3^2\vartheta_9^2+\vartheta_2^2\vartheta_{10}^2=0,\cr&
\vartheta_5^2\vartheta_7^2-\vartheta_1^2\vartheta_9^2+\vartheta_4^2\vartheta_{10}^2=0,\cr&
\vartheta_4^2\vartheta_7^2-\vartheta_3^2\vartheta_8^2-\vartheta_5^2\vartheta_{10}^2=0,\cr&
\vartheta_3^2\vartheta_7^2-\vartheta_4^2\vartheta_8^2-\vartheta_6^2\vartheta_9^2=0,\cr&
\vartheta_2^2\vartheta_7^2-\vartheta_1^2\vartheta_8^2-\vartheta_6^2\vartheta_{10}^2=0,\cr&
\vartheta_1^2\vartheta_7^2-\vartheta_2^2\vartheta_8^2-\vartheta_5^2\vartheta_9^2=0,\cr&
\vartheta_5^4-\vartheta_6^4-\vartheta_9^4+\vartheta_{10}^4=0,\cr}$$}
\vbox{\hsize=6cm
$$\leqalignno{&
\vartheta_4^2\vartheta_5^2-\vartheta_2^2\vartheta_6^2-\vartheta_7^2\vartheta_{10}^2=0,\cr&
\vartheta_3^2\vartheta_5^2-\vartheta_1^2\vartheta_6^2-\vartheta_8^2\vartheta_{10}^2=0,\cr&
\vartheta_2^2\vartheta_5^2-\vartheta_4^2\vartheta_6^2-\vartheta_8^2\vartheta_9^2=0,\cr&
\vartheta_1^2\vartheta_5^2-\vartheta_3^2\vartheta_6^2-\vartheta_7^2\vartheta_9^2=0,\cr&
\vartheta_3^4-\vartheta_4^4-\vartheta_6^4+\vartheta_{10}^4=0,\cr&
\vartheta_2^2\vartheta_3^2-\vartheta_1^2\vartheta_4^2+\vartheta_9^2\vartheta_{10}^2=0,\cr&
\vartheta_1^2\vartheta_3^2-\vartheta_2^2\vartheta_4^2-\vartheta_5^2\vartheta_6^2=0,\cr&
\vartheta_2^4-\vartheta_4^4-\vartheta_8^4+\vartheta_{10}^4=0,\cr&
\vartheta_1^2\vartheta_2^2-\vartheta_3^2\vartheta_4^2-\vartheta_7^2\vartheta_8^2=0,\cr&
\vartheta_1^4-\vartheta_2^4-\vartheta_6^4-\vartheta_9^4=0.\cr
}$$}}
\smallni
We order the 6 odd characteristics as follows:
$$(n^{(1)},\dots,n^{(6)})=\pmatrix{
0&0&1&1&1&1\cr
1&1&0&0&1&1\cr
0&1&1&1&0&1\cr
1&1&0&1&1&0\cr
}.$$
We use the notation
$$e(m)=(-1)^{\transpose a b}\qquad \bigl(m={a\choose b}\bigr).$$
Recall that a triplet $m_1,m_2,m_3$ of characteristics is called {\it azygetic\/} if
they are pairwise different and if
$$e(m_1)e(m_2)e(m_3)e(m_1+m_2+m_3)=-1.$$
The following result has been stated without  proof  by Rosenhain and 
proved  by Thomae and Weber. A proof can be found in [Fi].
\proclaim
{Lemma}
{For two different odd characteristics $m,n$ there exist 4 even characteristics
$n_1,\dots,n_4$ such that $m,n,n_i$ is azygetic. If we consider
the pair $(\grad\vartheta[m],\grad\vartheta[n])$ as a $2\times 2$ matrix then
$$\det(\grad\vartheta[m],\grad\vartheta[n])=\pm\pi^2\vartheta[n_1]\cdots\vartheta[n_4].$$
}
AleS%
\finishproclaim
Since the signs are essential for us, we collect them in a table. This table can be found in
[Fi]. (One sign in [Fi] had to be corrected.)
We use the abbreviation
$$D(i,j)=\det(\grad\vartheta[n^{(i)}],\grad\vartheta[n^{(j)}]),\quad 1\le i<j\le 6.$$
\halign{\qquad\qquad$#$,\hfil\quad&$#$,\hfil\quad&$#$\hfil\quad\cr
D(1,2)=\pi^{-2}\vartheta_7\vartheta_8\vartheta_9\vartheta_{10}&
D(1,3)=\vartheta_2\vartheta_3\vartheta_5\vartheta_7&
D(1,4)=\vartheta_1\vartheta_4\vartheta_5\vartheta_8,\cr
D(1,5)=-\vartheta_3\vartheta_4\vartheta_6\vartheta_{10}&
D(1,6)=\vartheta_1\vartheta_2\vartheta_6\vartheta_9&
D(2,3)=\vartheta_1\vartheta_4\vartheta_6\vartheta_7,\cr
D(2,4)=\vartheta_2\vartheta_3\vartheta_6\vartheta_8&
D(2,5)=-\vartheta_1\vartheta_2\vartheta_5\vartheta_{10}&
D(2,6)=\vartheta_3\vartheta_4\vartheta_5\vartheta_9,\cr
D(3,4)=-\vartheta_5\vartheta_6\vartheta_9\vartheta_{10}&
D(3,5)=-\vartheta_1\vartheta_3\vartheta_8\vartheta_9&
D(3,6)=\vartheta_2\vartheta_4\vartheta_8\vartheta_{10},\cr
D(4,5)=-\vartheta_2\vartheta_4\vartheta_7\vartheta_9&
D(4,6)=\vartheta_1\vartheta_3\vartheta_7\vartheta_{10}&
D(5,6)=\vartheta_5\vartheta_6\vartheta_7\vartheta_8.\cr
}
\medni
We describe 20 relations between the 6 generators $\grad_i$ of the module $\calM$.
We use the notation
$$\grad_i=\grad\vartheta[n_i],\quad 1\le i\le 6.$$
\proclaim
{Lemma}
{For $1\le i<j<k\le 6$ the relation
$$D(i,j)\grad_k=D(i,k)\grad_j-D(j,k)\grad_i$$
holds. Each of them is divisible by one of the $\vartheta_i$. Hence we obtain
$20$ relations where a typical one is
$$\vartheta_1\vartheta_4\vartheta_6\grad_1-\vartheta_2\vartheta_3\vartheta_5\grad_2-
\vartheta_8\vartheta_9\vartheta_{10}\grad_3=0.$$
}
RelD%
\finishproclaim
The proof is trivial. Just notice that the occurring $D$-s are just the
$2\times 2$ sub-determinants of the $2\times 3$-matrix $(\grad_i,\grad_j,\grad_k)$.
Now Lemma \RelD\ is just a consequence of the known fact that the cross product 
$a\times b$ of two
vectors in $\cz^3$ is orthogonal to both $a,b$.
\qed
\smallskip
We fix two odd characteristics $m,n$.
Every homogenous element of $\calM^+$ can be written in the form
$$T=g_m\grad\vartheta[m]+g_n\grad\vartheta[n]$$
with two meromorphic modular forms from $A(\Gamma[4,8])$. From Lemma
\AleS\ we can deduce that $g_i\vartheta[n_1]\cdots\vartheta[n_4]$ are holomorphic.
Hence $\calM$ is contained in
$$\calM(m,n):={1\over \vartheta[n_1]\cdots\vartheta[n_4]}(A(\Gamma[4,8])\grad\vartheta[m]+
A(\Gamma[4,8])\grad\vartheta[n]).$$
We can vary $m,n$ and obtain
$$\calM\subset\bigcap_{m,n}\calM(m,n).$$
The elements in the right hand side have poles outside the zeros of the forms
$\vartheta[n_1]\cdots\vartheta[n_4]$. Since these 15 forms have no joint zero in $\hz_2$, the
elements of the intersection
are holomorphic. Hence we obtain the following proposition.
\proclaim
{Proposition}
{We have
$$\calM=\bigcap_{m,n}\calM(m,n).$$}
IntrS%
\finishproclaim
We consider the submodule $\calN$ of $\calM$ that is generated by all
$\grad_i$. Proposition \IntrS\ shows that $\calM$ is a submodule of $(1/\chi_5)\calN$.
It is described as finite intersection of certain submodules which are defined 
by means of finitely many generators. As soon as we understand the structure
of $(1/\chi_5)\calN$, or equivalently, of $\calN$, we have a chance to determine
this intersection. In the next section we shall describe all relations between the
six $\grad_i$.
\neupara{Relations}%
In Lemma \RelD\ we described some of the relations between the $\grad_i$. 
It will turn all they do not generate all relations.
To describe all relations we introduce the free module of rank 6 over $A(\Gamma[4,8])$.
We denote the generators by $T_1,\dots,T_6$. We have to describe the
kernel of the natural homomorphism
$$\calF\lo\calN,\quad T_i\loma \grad_i.$$
We denote by $\calK$ the submodule of $\calF$ that is generated by the 20
elements which arise in Lemma \RelD. A typical example is
$$\vartheta_1\vartheta_4\vartheta_6T_1-\vartheta_2\vartheta_3\vartheta_5T_2-
\vartheta_8\vartheta_9\vartheta_{10}T_3$$
\proclaim
{Lemma}
{Let $T$ be an element of the kernel of $\calF\to\calN$. Then $\chi_5T$ is contained in
$\calK$. Hence the kernel of $\calF\to\calN$
equals the kernel of
$$\calF\buildrel{\cdot\chi_5}\over\lo\calF\lo\calF/\calK.$$}
RelE%
\finishproclaim
{\it Proof.\/} Let $P_1\grad_1+\cdots+P_6\grad_6=0$ be a (homogenous) relation.
After multiplication  by $\chi_5$
we can use the relations in Lemma \RelD\ to eliminate in this relation all $\grad_i$, $i>2$.
Then we obtain a relation between $\grad_1,\grad_2$. But these two forms
are independent due to Lemma \AleS. Hence the above relation, after multiplication by
$\chi_5$,  is a consequence
of the relations in Lemma \RelD.\qed
\smallskip
In principle, Lemma \RelE\ is a complete description of the module $\calN$. We can use it
to work out a finite  generating system of relations. For this,
we describe some extra relations between the $\grad_i$.
\proclaim
{Lemma}
{For each ordered pair of two different odd characteristics there is a relation
which is determined by this pair up to the sign. The relation
that belongs to pair $(5,6)$ in our numbering is
$$\vartheta_6^2D(1,5)\grad_1-\vartheta_5^2D(2,5)\grad_2-
\vartheta_8^2D(3,5)\grad_3+\vartheta_7^2D(4,5)\grad_4=0.$$
The full modular group acts
transitively on these $30$ relations (counted up to the sign).
In general the relation for a pair $(\alpha,\beta)$ is
the sum of four $\pm \vartheta_k^2D(i,\alpha)\grad_i$, $i\ne\alpha,\beta$,
where $\vartheta_k$ is the only theta that divides $D(i,\alpha)$ and
$D(\alpha,\beta)$.
}
ExtrA%
\finishproclaim
{\it Proof.\/} Along the lines of the proof of Lemma \RelE,
One multiplies the claimed relation by $\chi_5$ and eliminates $\grad_i$, $i>2$.
Then one obtains an expression $P_1\grad_1+P_2\grad_2$ with explicitly given
polynomials in the $\vartheta_i$. One has to show $P_1=P_2=0$.
We omit the straightforward calculation and mention only that for this one has
to use the Riemann relations.\qed
\smallskip
There is a second kind of extra relations between the $\grad_i$. To explain them,
we need some facts about theta characteristics. In [GS] it has been proved that
each odd characteristic $n$ can written in 12 different ways
(up to ordering) as a sum of five pairwise different even
characteristics
$$n=m_1+\cdots+m_5.$$
The full modular
group acts transitively on the set of all $\{n,m_1,\dots,m_5\}$.
For each of them we define the modular form
$$S:=S(n,m_1,\dots,m_5)=\vartheta[m_1]\cdots\vartheta[m_5]\grad\vartheta[n].$$
Hence we obtain 72 modular forms.
\smallskip
As we mentioned, the forms $\vartheta[m]$ and $\grad\vartheta[n]$
are modular forms with respect to the group $\Gamma_2[2,4]$.
As a consequence, we get
$$S\in[\Gamma_2[2,4],3,\chi_S,\St]$$
with a certain quadratic character $\chi_S$ on $\Gamma_2[2,4]$.
The information about these characters which we need can be taken
from the paper [SM].
\proclaim
{Lemma}
{The $72$ forms $S$ are modular forms with respect to the group
$\Gamma_2[2,4]$ and a certain quadratic character $\chi_S$.
In this  way there arise $12$ different characters and to the 
ach associated space of  modular 
forms belong  six of the  forms  S.
Each of the $72$ forms $S$ is uniquely determined by its odd characteristic
$n$ and the character $\chi_S$.}
CharChar%
\finishproclaim
We denote the form $S$ which belongs to $n$ and $\chi$ by $S(n,\chi)$.
\smallskip

\proclaim
{Lemma}
{
We fix one of the $12$ characters $\chi$. Let $S_1,\dots,S_6$ be the
six functions $S(n,\chi)$. If one cancels one of the six, say $S_6$,
one gets a relation between the other five of the following type.
$$\sum_{i=1}^5\pm\vartheta[m_i]^2S_i=0.$$
Here $m_i$ are certain even characteristics which are uniquely determined
and also the signs (up to a common sign change) are uniquely determined.
}
ExtrB%
\finishproclaim
The rule how the $m_i$ can be found is a little complicated.
We explain how $m_1$ can be found.
\smallni
1) There are three even characteristics $q_1,q_2,q_3$ that occur in
$S_1$ but not in $S_6$ (the form which has been cancelled).
\vskip1mm\noindent
2) There is one pair in $\{q_1,q_2,q_3\}$, say $\{q_2,q_3\}$, such that
$\vartheta[q_2]\vartheta[q_3]$ does not occur in any of the four $S_2,\dots,S_5$.
\vskip1mm\noindent
Then one has to use $m_1:=q_1$ in Lemma \ExtrB.
\smallskip
We explain this in an example.
The six forms
$$\eqalign{
&\hbox to 4cm{$
S_1=\vartheta_3\vartheta_5\vartheta_6\vartheta_8\vartheta_9\grad_1$,\hfil}
\hbox to 4cm{$
S_2=\vartheta_1\vartheta_2\vartheta_4\vartheta_8\vartheta_9\grad_2$,\hfil}
\hbox{$
S_3=\vartheta_2\vartheta_5\vartheta_7\vartheta_8\vartheta_{10}\grad_4$,\hfil}\cr
&\hbox to 4cm{$
S_4=\vartheta_4\vartheta_6\vartheta_7\vartheta_9\vartheta_{10}\grad_6$,\hfil}
\hbox to 4cm{$
S_5=\vartheta_1\vartheta_3\vartheta_4\vartheta_5\vartheta_{10}\grad_3$,\hfil}
\hbox{$
S_6=\vartheta_1\vartheta_2\vartheta_3\vartheta_6\vartheta_7\grad_5$\hfil}
\cr}
$$
belong to the same character. The characteristic that occur in $S_1$ but not
in $S_6$ are $m_5,m_8,m_9$. The pair $m_5,m_9$ does not occur in $S_2,\dots,S_5$.
Hence the theta square which has to be added is $\vartheta[m_8]^2$.
The whole relation will be
$$\vartheta[m_8]^2S_1-\vartheta[m_9]^2S_2+\vartheta[m_5]^2S_3-\vartheta[m_4]^2S_4+
\vartheta[m_{10}]^2S_5.$$
The determination of the signs needs some extra work.
\smallskip
We do not give details of the proof of Lemma \ExtrB. We only mention that
it is similar to the proof of Lemma \ExtrA.
\smallskip
The relations that we described so far are defining relations.
\proclaim
{Proposition}
{The module of
relations between the six forms $\grad_i$, i.e.~the kernel of the natural
homomorphism $\calF\to\calN$, is generated by the $20$ relations described
in Lemma \RelD, the $30$ relations described in Lemma \ExtrA\ and the $72$
relations described in Lemma \ExtrB.}
AllRel%
\finishproclaim
Using Lemma \RelE, the proof can be given with the help of a computer.\qed
\neupara{A structure theorem}%
Now we have the possibility to determine the structure of $\calM$. From Proposition \IntrS\
we get
$$\chi_5\calM=\bigcap_{m,n}
{\chi_5\over \vartheta[n_1]\cdots\vartheta[n_4]}(A(\Gamma[4,8])\grad\vartheta[m]+
A(\Gamma[4,8])\grad\vartheta[n]).
$$
The right-hand side is a submodule of $\calN$ which we understand completely
(Proposition \AllRel).
Hence it is possible to compute the intersection with the help of a computer.
We did this by means of the computer algebra system {\pro SINGULAR}. In this way
we could determine a finite system of generators of $\calM$ and we also could get
the Hilbert function. We mention that $\calM$ is bigger than $\calN$.
We have to describe now the extra generators.
\proclaim
{Proposition}
{The modular form
$$
{(\vartheta_4\vartheta_6^4
\vartheta_8
+\vartheta_4
\vartheta_8\vartheta_9^4)\grad_1
-\vartheta_1\vartheta_6
\vartheta_9\vartheta_{10}^3\grad_3\over\vartheta_2\vartheta_5}$$
is holomorphic, hence contained in $[\Gamma[4,8],2,\St]$.
It is not contained in the submodule $\calN$. The orbit under the full
modular form consists up to constant factors of $360$ modular forms.
}
ExtrH%
\finishproclaim
By means of {\pro SINGULAR} one can verify that the $A(\Gamma[4,8])$-module
generated by the six
$\grad_i$ and
the $\Gamma_2$-orbit of the form described in Proposition \ExtrH\ equals the
module $\calM$. {\pro SINGULAR} also gives the Hilbert function.
\proclaim
{Theorem}
{The $A(\Gamma[4,8])$-module
$$\calM=\bigoplus [\Gamma[4,8], r/2 , v_{\vartheta}^r, \St]$$
is generated by the six
$\grad_i$ and
the $\Gamma_2$-orbit of the form described in Proposition \ExtrH.
The Hilbert function is
$$\eqalign{&
\quad\sum_{r=0}^\infty [\Gamma[4,8], r/2 , v_{\vartheta}^r, \St]t^r=\cr&
{60t^9 - 60t^8 - 318t^7 + 252t^6 + 606t^5 + 316t^4 + 126t^3 + 36t^2 +
6t\over(1-t)^4}=\cr&
6t + 60t^2 + 330t^3 + 1300t^4 + 4060t^5 + 9952t^6 + 20000t^7 +
    35168t^8 +\cdots\cr
}$$
}
MainT%
\finishproclaim
\vskip1.5cm\noindent
{\paragratit References}
\bigskip
\item{[Ao]} Aoki, H.:
{\it On Vector Valued {S}iegel modular forms of degree 2 with small levels,\/}
Osaka Journal of Mathematics {\bf 49}, 625--651 (2012)
\medskip
\item{[CG]} Cl{\'e}ry, F., van der Geer,~G., Grushevsky,~S.:
{\it Siegel modular forms of genus 2 and level 2,\/}
math.AG/1306.6018  (2013)
\medskip
\item{[Do]} von Dorp,~C.:
{\it Vector-valued Siegel modular forms of genus 2,\/}
MSc Thesis, Korteweg-de Vries Instituut voor Wiskunde,
Universiteit van Amsterdam (2011)
\medskip
\item{[Fi]} Fiorentino, A.:
{\it On a ring of modular forms related to the Theta gradients map in genus 2,\/}
Thesis at the Universit\`a di Roma la Sapienza, Italy,\hfill\break
arXiv: 1109.2362 (2011)
\medskip
\item{[Fr]} Freitag, E.: {\it Siegelsche Modulfunktionen,} Grundlehren
der mathematischen Wissenschaften, Bd. {\bf 254}. Berlin Heidelberg New
York: Springer (1983)
\medskip
\item{[GG]} van Geemen,~B., van der Geer,~G.: {\it Kummer  varieties and the
moduli spaces of abelian varieties,\/}, Amer. J. of Math. {\bf 108}, 615--642 (1986)
\medni
\item{[GS]} van Geemen,~B., van Straten,~D.:
{\it The cusp forms of weight  3 on $\Gamma_2(2,4,8)$.\/}
Mathematics of Computation, {\bf 61}, Number 204, 849--872  (1993)
\medskip
\item{[Ib]} Ibukiyama,~T.:
{\it Vector valued Siegel modular forms of symmetric tensor weight of small degrees,\/}
Comment.\ Math.\ St.\ Pauli {\bf 61}, No.~1, 51--75 (2012)
\medskip
\item{[Ig]} Igusa, J.I.:
{\it On Siegel modular forms of genus II.\/}
Amer.\ J.\ Math.\ {\bf 86}, 392--412 (1964)
\medskip
\item{[Ru1]} Runge, B.:
{\it On Siegel modular form, part I\/}, J. Reine angew. Math. {\bf 436}, 57-85 (1993)
\medskip
\item{[Ru2]} Runge, B.:
{\it On Siegel modular forms, part II\/}, Nagoya Math. J. {\bf138}, 179-197 (1995)
\medskip
\item{[SM]} Salvati Manni, R.:
{\it Modular varieties with level 2 theta structure.\/}
Am.\ Journal pf Math.\ {\bf 116}, No 6, 1489--1511 (1994)
\medskip
\item{[Sa]} Sasaki,~R.: {\it Modular forms vanishing at the reducible
points of the Siegel upper-half space,\/} J. Reine Angew. Math. {\bf 345}, 111--121 (1983).
\medskip
\item{[Sat]} Satoh,~T.: {\it On Certain Vector Valued Siegel Modular forms
of Degree Two,\/}  Mathematische
Annalen {\bf274}, 365--387 (1986)
\medskip
\item{[Wi]} Wieber,~T.:
{\it Structure Theorems for Certain Vector Valued {S}iegel Modular Forms of Degree Two,\/}
doctoral thesis, University of Heidelberg  (2013)
\bye